\documentclass[12pt]{amsart}
\usepackage{a4,amscd,amssymb,amsfonts,verbatim}

\newcommand{\Co}{{\mathbb C}}
\renewcommand{\Re}{{\mathbb R}}
\newcommand{\Za}{{\mathbb Z}}

\newcommand{\n}{\nabla}      
\newcommand{\tn}{\widetilde{\nabla}}

\newcommand{\bD}{\overline{D}}

\newcommand{\Si}{\operatorname{Spin}}

\renewcommand{\phi}{\varphi}
\newcommand{\one}{{\bf 1}}

\newcommand{\tR}{\widetilde{R}}

\newcommand{\half}{\frac{1}{2}}

\newcommand{\hA}{\widehat{A}}

\newcommand{\tM}{\widetilde{M}}
\newcommand{\tl}{\tilde{l}}
\newcommand{\bM}{\overline{M}}
\newcommand{\bareta}{\overline{\eta}}

\newcommand{\bS}{\overline{S}}
\newcommand{\grad}{\operatorname{grad}}
\newcommand{\spec}{\operatorname{spec}}
\newcommand{\sign}{\operatorname{sign}}
\newcommand{\Sign}{\operatorname{Sign}}
\newcommand{\chern}{\operatorname{ch}}
\newcommand{\trace}{\operatorname{tr}}

\newcommand{\la}{\langle}       
\newcommand{\ra}{\rangle}

\newtheorem{thm}{Theorem}[section]

\newtheorem{prop}[thm]{Proposition}
\newtheorem{cor}[thm]{Corollary}

\begin{document}

\title[dependence on the spin structure of the eta and rokhlin invariants]
{dependence on the spin structure of the eta and rokhlin invariants}

\author{Mattias Dahl}

\address{
Department of Mathematics \\
Royal Institute of Technology \\
S-100 44 Stockholm, Sweden\\
}

\email{dahl@math.kth.se}


\date{\today}

\begin{abstract}
We study the dependence of the eta invariant $\eta_D$ on the spin structure,
where $D$ is a twisted Dirac operator on a $(4k+3)$-dimensional spin manifold. 
The difference between the eta invariants for two spin structures related by
a cohomolgy class which is the reduction of a $H^1(M,\Za)$-class is
shown to be a half integer. 
As an application of the technique of proof the generalized Rokhlin invariant 
is shown to be equal modulo $8$ for two spin structures related in this way. 
\end{abstract}

\maketitle

\section{Introduction}

\subsection{The eta invariant}

Let $M$ be a compact manifold of dimension $4k+3$.
Let $D$ be a self-adjoint first order elliptic operator on $M$.
The eta invariant of $D$ is defined as 
$$
\eta_D = \half \Big( \eta(0) + \dim \ker D \Big)
$$
where the eta-function $\eta(s)$ is given by 
$$
\eta(s) = \sum_{\lambda \in \spec D, \lambda \neq 0} 
\frac{\sign \lambda}{|\lambda|^s}
$$
which converges for $\operatorname{Re} s$ large, and has an analytic
continuation which is regular at $s=0$, see \cite{Gilkey}.

Assume that $M$ is a spin manifold and let $\sigma$ be a spin
structure on $M$ with associated spinor bundle $S_{\sigma}$. 
Let $E$ be a vector bundle with connection on $M$. 
We use the notation $\eta(\sigma;E)$ for the eta invariant of 
the twisted Dirac operator $D^E$ on $S_{\sigma} \otimes E$. 
The set $\Si(M)$ of spin structures on $M$ is an affine space 
modeled on $H^1(M;\Za_2)$, i.e. the vector space $H^1(M;\Za_2)$ 
acts freely and transitively on $\Si(M)$:
$$
\Si(M) \times H^1(M;\Za_2) \ni (\sigma, \delta) \mapsto 
\sigma + \delta \in \Si(M).
$$
For $\sigma \in \Si(M)$ and $\delta, \delta' \in H^1(M;\Za_2)$
we define the first and second difference functions of $\eta(\sigma;E)$ by
$$
\Delta \eta(\sigma, \delta;E) = 
\eta(\sigma + \delta;E) - \eta(\sigma;E)
$$
and
$$
\Delta^2 \eta(\sigma,\delta',\delta;E) = 
\Delta\eta(\sigma + \delta',\delta;E) - \Delta\eta(\sigma,\delta;E).
$$
The first difference $\Delta \eta$ is a special case of the relative
eta invariant introduced by Atiyah, Patodi and Singer in \cite{APS2}
and since studied by many authors.
The ideas behind the following 
theorem also goes back to Atiyah, Patodi and Singer.
\begin{thm} \label{eta}
Let $(M,g)$ be a compact spin manifold of dimension $4k+3$.
Let $\sigma$ be a spin structure on $M$ and let 
$\delta, \delta' \in H^1(M;\Za_2)$.
Suppose $\delta$ is the reduction modulo $2$ of an integer class.
Let $E$ be a vector bundle on $M$ with connection $\n^E$.
Then 
\begin{enumerate}

\item[(a)] $\Delta \eta(\sigma,\delta;E) \in \half \Za$, 

\item[(b)] $\Delta^2 \eta(\sigma,\delta',\delta;E) \in \Za$.

\end{enumerate}

These differences do not depend on the metric $g$.
\end{thm}
 
\subsection{The Rokhlin invariant}

Let $M$ be a compact spin manifold of dimension $8k + 3, k \geq 0$, 
and suppose that $M$ is a boundary. The Rokhlin invariant is 
a function defined on the set $\Si(M)$ of spin structures on $M$ and 
taking values in the integers modulo 16,
$$
R_M: \Si(M) \to \Za_{16},
$$
see \cite{Turaev} and \cite{Finashin}.
Given a spin structure $\sigma$ on $M$, we take a spin manifold $N$
with spin structure $\tau$ such that $(M,\sigma)$ is the boundary 
of $(N,\tau)$. Such an $N$ exists since the forgetful homomorphism
of cobordism rings $\Omega^{\Si}_{8k+3} \to \Omega^{SO}_{8k+3}$ is
injective, see \cite[p.351]{Stong}, and thus $M$ is an oriented 
boundary if and only if it is a spin boundary. The Rokhlin function
is defined by 
$$
R_M(\sigma) = \Sign(N) \mod 16
$$
where $\Sign(N)$ is the signature of $N$. This is independent of
the choice of $N$ since the signature invariant is additive and
by a theorem of Ochanine, see \cite[p.113]{Hirzebruch-et.al}, 
the signature of a closed $8k+4$ dimensional spin manifold is 
a multiple of $16$. 

Lee and Miller \cite{Miller/Lee} showed that if we equip $M$ with 
a riemannian metric then the Rokhlin function can be computed in 
terms of eta invariants of the signature operator and of twisted 
Dirac operators,
\begin{equation} \label{roketa}
R_M(\sigma) = - \eta_{\operatorname{Hirz}} +
8 \sum_i b_i \eta_{D^{Z_i}} \mod 16,
\end{equation}
where $b_i$ are integers and $Z_i$ are vector bundles.
In dimension $3$ this formula is
\begin{equation} \label{dim3}
R_M(\sigma) = - \eta_{\operatorname{Hirz}} 
- 8 \eta_D \mod 16,
\end{equation}
and in dimension $11$ we have 
\begin{equation} \label{dim11}
R_M(\sigma) = - \eta_{\operatorname{Hirz}} 
+ 8 \eta_{D^{TM}} - 32 \eta_D \mod 16.
\end{equation}
Since the expression on the right-hand side of (\ref{roketa}) is
defined even if $M$ is not a boundary, we take (\ref{roketa}) as the 
definition of $R_M$ for a general $8k+3$-dimensional spin manifold. 
From (\ref{roketa}) it also follows that this is independent of the 
riemannian metric, but a priori this extension of the Rokhlin
function takes it values in $\Re$ modulo $16$. A theorem by Fischer and 
Kreck \cite{Fischer/Kreck} tells us that (\ref{roketa}) is always an 
integer modulo $16$.

In this paper we are going to use (\ref{roketa}) to study how the 
Rokhlin invariant varies with the spin structure on $M$. 
For $\sigma \in \Si(M)$ and $\delta_0, \dots,\delta_m \in H^1(M;\Za_2)$ 
the difference functions of $R_M$ are defined inductively by
$$
\begin{aligned}
\Delta^1 R_M(\sigma, \delta_0) =& R_M(\sigma + \delta_0) -
R_M(\sigma), \\
\Delta^{m+1} R_M(\sigma, \delta_0,\dots,\delta_m) =& 
\Delta^m R_M(\sigma + \delta_0,\delta_1,\dots,\delta_m) \\
&- \Delta^m R_M(\sigma, \delta_1,\dots,\delta_m).
\end{aligned}
$$
In \cite{Turaev} Turaev shows that $\Delta^4 R_M = 0$ for any 
3-dimensional manifold $M$. Finashin conjectures in \cite{Finashin} 
that $\Delta^4 R_M = 0$ holds in all dimensions.

Our main result is the following.

\begin{thm} \label{main}
If $\delta \in H^1(M;\Za_2)$ is the reduction modulo 2 of an 
integer class then
$\Delta R_M(\sigma,\delta) = 0 \text{ or } 8 \mod 16$.
\end{thm}

This generalizes a result by Taylor, 
see \cite[Theorem 6.2.]{Taylor}, in dimension 3. In the case where 
$M$ is a boundary this follows by results of Finashin in
\cite{Finashin} but in the general case it seems to be new.

\section{Spin structures and Dirac operators}

Throughout this paper $(M,g)$ will denote a compact oriented
riemannian spin manifold of dimension $n$. For the necessary 
background on spin geometry we refer to \cite{Gilkey} and \cite{LM}.

\subsection{Spin structures and spinor bundles}

Let $\pi: \Si(n) \to SO(n)$ be the double cover of the group $SO(n)$.
A spin structure $\sigma$ on $M$ is a $\Si(n)$-principal bundle 
over $M$ which is pointwise a double cover of the oriented orthonormal 
frame bundle $SO(M)$,
$$
\xi : \sigma \to SO(M), 
$$
$$
\xi(pg) = \xi(p)\pi(g), \quad p \in \sigma, g \in \Si(n).  
$$
We denote the set of (isomorphism classes of) spin structures 
by $\Si(M)$. This is an affine space modeled on $H^1(M;\Za_2)$, 
given two spin structures $\sigma, \sigma'$ there is an element 
$\sigma' - \sigma = \delta \in  H^1(M;\Za_2)$ which is the 
difference of the spin structures. This difference element can 
be constructed as follows. Let $\{ U_i, a_{ij} \}$ be a
trivialization of the bundle $SO(M)$ with transition functions
$$
a_{ij} : U_i \cap U_j \to SO(n)
$$
and let 
$$
b_{ij}, b'_{ij} : U_i \cap U_j \to \Si(n)  
$$
be the transition functions for $\sigma, \sigma'$.
Then 
$$
\pi \circ b_{ij} = \pi \circ b'_{ij} = a_{ij}
$$
and the $c_{ij}$ defined as 
\begin{equation} \label{trans-rel}
c_{ij} = b'_{ij} b^{-1}_{ij} : U_i \cap U_j \to \ker(\pi) 
= \{ \pm 1 \} = \Za_2
\end{equation}
give $\delta \in  H^1(M;\Za_2)$. The element 
$\delta$ can be viewed as a $\Za_2$-principal fibre bundle.
By the action of $\Za_2 = \{ \pm 1 \}$ on $\Co$ we get an 
associated flat complex line bundle which we denote by $L_{\delta}$.

Let $\rho : \Si(n) \to \operatorname{End}({\mathcal S})$ be the
spinor representation of the spin group on the vector space ${\mathcal S}$.
We denote by $S_{\sigma}$ the spinor bundle 
$\sigma \times_{\rho}{\mathcal S}$ associated to $\sigma$. 
In the trivialization introduced above $S_{\sigma}$ will have
transition functions 
$\rho(b_{ij}) : U_i \cap U_j \to \operatorname{End}({\mathcal S})$
and from (\ref{trans-rel}) it follows that the transition functions 
for $S_{\sigma'}$ are
$$
\rho(b'_{ij}) = \rho(c_{ij})\rho(b_{ij}).
$$
Since $\rho(-1) = -1$ this translates to
\begin{equation} \label{spin-repr}
S_{\sigma'} = L_{\delta} \otimes S_{\sigma}.
\end{equation}

In this paper we will be interested in the case when the line 
bundle $L_{\delta}$ is topologically trivial, that is when 
$L_{\delta}$ has a non-vanishing section. A smooth non-vanishing 
section $l$ of $L_{\delta}$ defines an invertible map
$$
l: \phi \mapsto l \otimes \phi 
$$
which sends sections of $S_{\sigma}$ to sections of $S_{\sigma'}$. 
If we normalize $l$ to $|l|=1$ this map is pointwise an isometry since
$$
\la l\otimes \phi, l \otimes \psi \ra = |l|^2 \la \phi, \psi \ra
$$
and also a surjective isometry of the Hilbert spaces of square 
integrable sections
$$
l: L^2(S_{\sigma}) \to L^2(S_{\sigma'}).
$$

Let $H^1(M;\Za_2)_0$ be the subspace of $H^1(M;\Za_2)$ 
consisting of $\delta$ such that $L_{\delta}$ is topologically 
trivial. It can be characterized as follows, 
see \cite[p.84]{LM} and \cite{Ammann}.
\begin{prop} \label{alphaprop}
$\delta \in H^1(M;\Za_2)_0$ if and only if there is an integral 
$\alpha \in H^1(M;\Re)$ such that
$$
\delta([\gamma]) = e^{\pi i \int_{\gamma} \alpha}
$$
for all closed curves $\gamma$. 
\end{prop}

\begin{proof}
Suppose $l$ is a section of $L_{\delta}$ with $|l|=1$. 
Such a section is the same as a $\delta$-equivariant
function $\tl = e^{\pi i \theta}$ on the universal cover 
$\tM$ of $M$. Let $\alpha = d\theta = \frac{1}{\pi i} l^{-1} dl$.
Then $\alpha$ drops down to a closed real-valued one-form on $M$
and
\begin{eqnarray*}
e^{\pi i \int_{\gamma} \alpha} &=& 
e^{\pi i \int_{\tilde{\gamma}} d\theta} \\ 
&=& e^{\pi i \Big( \theta(\gamma(x_0)) - \theta(x_0) \Big)} \\
&=& \tl(\gamma(x_0)) \tl(x_0)^{-1} \\
&=& \delta([\gamma])
\end{eqnarray*}
where $\tilde{\gamma}$ is the lift of $\gamma$ to $\tM$, a path 
from $x_0$ to $\gamma(x_0)$. On the other hand, suppose 
$\delta([\gamma]) = e^{\pi i \int_{\gamma} \alpha}$
for some closed one-form $\alpha$. Then the pullback of 
$\alpha$ to $\tM$ is the differential of some function $\theta$
and $\tl = e^{\pi i \theta}$ is a $\delta$-equivariant
function on $\tM$, which gives a non-vanishing section of 
$L_{\delta}$. 
\end{proof}

\subsection{Dirac operators}
For a spin structure $\sigma$ with associated spinor bundle
$S_{\sigma}$ the Dirac operator $D$ is defined by 
$$
D\phi = e_i \cdot \n_{e_i} \phi 
$$
where $e_i$ is a local orthonormal frame, the dot is Clifford 
multiplication, and $\phi$ is a smooth section of $S_{\sigma}$. 
$D$ extends to a self-adjoint operator on $L^2(S_{\sigma})$.
As before we let $\sigma$ and $\sigma' = \sigma + \delta$ be 
spin structures on $M$ and we assume that $\delta \in H^1(M;\Za_2)_0$.
We let $l$ be a unit norm section of $L_{\delta}$ and define 
a connection $\n'$ by
\begin{equation} \label{nablaprim}
\n' \phi = l^{-1} \n (l\phi) = \n \phi + l^{-1} dl \phi
\end{equation}
acting on sections of $S_{\sigma}$. This connection is metric:
\begin{eqnarray*}
X \la \phi,\psi \ra &=&  X \la l\phi,l\psi \ra  \\
&=& \la \n_X (l\phi),l\psi \ra + \la l\phi,\n_X(l\psi) \ra \\
&=& \la l^{-1}\n_X (l\phi),\psi \ra + \la \phi,l^{-1}\n_X(l\psi) \ra \\
&=& \la \n'_X\phi,\psi \ra  + \la \phi,\n'_X\psi \ra.
\end{eqnarray*}
Let $D'$ be the operator of Dirac type defined using $\n'$,
$$
D' = e_i \cdot \n'_{e_i} = l^{-1} D l = D +  l^{-1} \grad l\cdot.
$$
The following is a fundamental observation. 
\begin{prop}
$D'$ acting on $L^2(S_{\sigma})$ is isospectral to 
$D$ acting on $L^2(S_{\sigma'})$. 
\end{prop}
\begin{proof}
Suppose $\phi$ is a section of $S_{\sigma'}$ and an
eigenspinor of $D$ with eigenvalue $\lambda$. Then 
$$
D'(l^{-1}\phi) = l^{-1} D \phi = \lambda l^{-1}\phi
$$ 
so $l^{-1}\phi$ is an eigenspinor of $D'$, also with eigenvalue
$\lambda$. If $\psi$ is an eigenspinor of $D'$ then $l\psi$ is 
an eigenspinor of $D$ with the same eigenvalue.
\end{proof}

This means that $D'$ acting on $L^2(S_{\sigma})$ has the same 
spectral invariants as $D$ acting on $L^2(S_{\sigma'})$. 

\section{The eta invariant}

In this section we prove the theorem on the eta invariant stated in
the introduction and state some corollaries.

\begin{thm} \label{eta1}
Let $(M,g)$ be a compact spin manifold of dimension $4k+3$.
Let $\sigma$ be a spin structure on $M$ and let $\delta \in H^1(M;\Za_2)_0$,
$\delta' \in H^1(M;\Za_2)$.
Let $E$ be a vector bundle on $M$ with connection $\n^E$.
Then 
\begin{enumerate}

\item[(a)] $\Delta \eta(\sigma,\delta;E) \in \half \Za$, 

\item[(b)] $\Delta^2 \eta(\sigma,\delta',\delta;E) \in \Za$.

\end{enumerate}

These differences do not depend on the metric $g$.
\end{thm}
 
\begin{proof}
As before let $l$ be a normalized section of $L_{\delta}$ and let
$\n' = l^{-1} \n l$. Let $\bM = M \times I$ where $I = [0,1]$
and denote by $\nu$ the vectorfield $\frac{\partial}{\partial t}$ 
on $\bM$.
The spin structure $\sigma$ on $M$ induces a spin structure on 
$\bM$, the associated spinor bundle is equal to 
$S_{\sigma} \oplus S_{\sigma}$. To shorten the notation we will write
this as $\bS$.

Let $\chi: I \to I$ be a smooth function such that $\chi(t)= 0, t \leq 1/3$ 
and $\chi(t)= 1, t \geq 2/3$. Define a connection $\tn$ 
on sections of the trivial complex line bundle $\one$ over $\bM$ by 
$$
\tn_X f = \Big( X + \chi(t) l^{-1}dl(X) \Big) f  , 
\quad
\tn_{\nu} f = \nu f
$$
where $X \in TM$. This connection is metric and the induced connection
$\n \otimes \tn$ on $\bS \otimes \one$ interpolates between $\n$ 
close to $t=0$ and $\n'$ close to $t=1$. 

Let $\bD^E$ be the Dirac operator on $\bS \otimes E \otimes \one$ 
constructed using the connection $\n \otimes \n^E \otimes \tn$.
Close to the boundary component $M \times \{ 0 \}$ of $\bM$ we have 
$\bD^E = \nu \cdot (\frac{\partial}{\partial t} + D^E)$ and close to 
$M \times \{ 1 \}$ we have 
$\bD^E = \nu \cdot (\frac{\partial}{\partial t} + D^{'E})$. 

The Atiyah-Patodi-Singer index theorem \cite[p.305]{Gilkey} tells us that 
\begin{equation}
\label{APS1}
\operatorname{ind}(\bD^E) = 
\int_{\bM} \hA(M) \wedge \chern(\n^E) \wedge \chern (\tn) 
- \eta_{D^E} + \eta_{D^{'E}}.
\end{equation}
where the index on the left hand side is the Fredholm index of $\bD^E$ 
acting on sections of $\bS \otimes \one \otimes E$ satisfying the 
Atiyah-Patodi-Singer boundary condition.

Since $D^{'E}$ acting on sections of $S_{\sigma} \otimes E$ is 
isospectral to $D^E$ acting on sections of $S_{\sigma + \delta} \otimes E$
we have $\eta_{D^{'E}} - \eta_{D^E} = \eta(\sigma+\delta;E) - 
\eta(\sigma;E) = \Delta\eta(\sigma,\delta;E)$ and
\begin{equation}
\label{APS2}
\Delta\eta(\sigma,\delta;E) =
- \operatorname{ind}(\bD^E) + 
\int_{\bM} \hA(M) \wedge \chern(\n^E) \wedge \chern(\tn). 
\end{equation}

Next we compute $\chern(\tn)$ and for this we need the curvature 
$\tR$ of $\tn$. For $X \in TM$ we write
$$
\tn_X = X + \chi(t) l^{-1} dl(X) = \n + \pi i \chi(t) \alpha(X)
$$
where as in the proof of Proposition (\ref{alphaprop}) 
$\alpha =  \frac{1}{\pi i} l^{-1} dl $ is a closed integer-valued
one-form. For $X,Y \in TM$ we have
\begin{eqnarray*}
\tn_X \tn_Y &=& 
\Big(X + \pi i \chi(t) \alpha(X)\Big)
\Big(Y + \pi i \chi(t) \alpha(Y)\Big) \\
&=& X Y + 
\pi i \chi(t)\Big(X(\alpha(Y)) + \alpha(Y) X + \alpha(X) Y\Big)  \\
&+& \Big(\pi i \chi(t)\Big)^2 \alpha(X)\alpha(Y)
\end{eqnarray*}
so 
\begin{eqnarray*}
\tR(X,Y) &=& \tn_X \tn_Y - \tn_Y \tn_X - \tn_{[X,Y]} \\
&=& \pi i \chi(t)
\Big(X(\alpha(Y)) - Y(\alpha(X)) - \alpha([X,Y]) \Big) \\
&=& \pi i \chi(t) d\alpha(X,Y) \\
&=& 0
\end{eqnarray*}
since $\alpha$ is closed. Next we have 
\begin{eqnarray*}
\tR(\nu,X) \phi &=& (\tn_{\nu} \tn_X - \tn_X \tn_{\nu})\phi \\
&=& \pi i \chi'(t) \alpha(X)
\end{eqnarray*}
and we conclude that
$$
\tR = \pi i d\chi \wedge \alpha.
$$
Thus we find that the Chern character of $\tn$ is
\begin{equation} \label{bchern}
\chern(\tn) = \trace e^{ \frac{i \tR}{2\pi}} = 
e^{- \frac{d\chi \wedge \alpha}{2}} =
1 - \frac{d\chi \wedge \alpha}{2} 
\end{equation} 
and we see that the integral in (\ref{APS1}) is 
$$
\int_{\bM} \hA(M) \wedge \chern(\n^E) \wedge \chern(\tn)
= 
\int_{\bM} \hA(M) \wedge \chern(\n^E) \wedge 
(1 - \frac{d\chi \wedge \alpha}{2}).
$$
When we integrate over the $I$ factor only the term with $d\chi$ will
contribute and we have left
\begin{equation} \label{APSintegral}
\int_{\bM} \hA(M) \wedge \chern(\n^E) \wedge \chern(\tn)
=
-\half \int_M \hA(M) \wedge \chern(\n^E) \wedge \alpha.
\end{equation}

The same calculation can now be used again to see that this integral 
is the index of a twisted Dirac operator. Let $S^1$ be the circle 
of length $1$. Then
\begin{eqnarray*}
\int_M \hA(M) \wedge \chern(\n^E) \wedge \alpha
&=&
\int_{M \times S^1}  \hA(M) \wedge \chern(\n^E) \wedge dt \wedge \alpha \\
&=& 
\int_{M \times S^1}  \hA(M) \wedge \chern(\n^E) \wedge (1+dt \wedge \alpha) \\
&=&
\int_{M \times S^1}  \hA(M) \wedge \chern(\n^E) \wedge e^{dt \wedge \alpha}.
\end{eqnarray*}
Since $dt \wedge \alpha \in H^2(M \times S^1;\Za)$
there is a complex line bundle $K \to M \times S^1$ with 
$c_1(K) = dt \wedge \alpha$, curvature 2-form 
$\frac{2 \pi}{i} dt \wedge \alpha$, and $\chern(\n^K) = e^{dt \wedge \alpha}$.
This means that (\ref{APSintegral}) is equal to 
\begin{equation} \label{APSintegral2}
\begin{aligned}
-\half \int_M \hA(M) \wedge \chern(\n^E) \wedge \alpha
&=
-\half \int_{M \times S^1}  \hA(M) \wedge \chern(\n^E) \wedge \chern(\n^K) \\
&= 
-\half \operatorname{ind}(D^{E \otimes K})
\end{aligned}
\end{equation}

From (\ref{APS2}) and (\ref{APSintegral2}) we now get 
\begin{equation} \label{deltaeta}
\Delta\eta(\sigma,\delta;E) =
- \operatorname{ind}(\bD^E) 
-\half \operatorname{ind}(D^{E \otimes K}),
\end{equation}
which is a half integer. We have proved part (a) of the theorem. 
To prove part (b) note that the last term,
$$
\operatorname{ind}(D^{E \otimes K}) = 
\int_M \hA(M) \wedge \chern(\n^E) \wedge \alpha,
$$
in (\ref{deltaeta}) does not depend on the spin structure $\sigma$. 
When we take a second difference,
$$
\Delta^2 \eta(\sigma,\delta',\delta;E) = 
\Delta\eta(\sigma + \delta',\delta;E) - \Delta\eta(\sigma,\delta;E),
$$
these terms will cancel, and part (b) follows.
\end{proof} 

In the case where the twisting bundle $E$ is flat, for example 
if $E$ is the trivial line bundle, we get the following corollary.
\begin{cor}
Let $(M,g)$ be a compact spin manifold of dimension $4k+3$.
Let $\sigma$ be a spin structure on $M$ and suppose 
$\delta \in H^1(M;\Za_2)_0$. 
Let $E$ be a vector bundle on $M$ with a flat connection $\n^E$.
Then $\Delta \eta(\sigma,\delta;E) \in \Za$,
\end{cor}
\begin{proof}
Since $\chern(\n^E) = 1$ we have from (\ref{APS2}) and
(\ref{APSintegral}) 
\begin{eqnarray*}
\Delta \eta(\sigma,\delta;E) 
&=&  
- \operatorname{ind}(\bD) 
-\half \int_M \hA(M) \wedge \chern(\n^E) \wedge \alpha \\
&=& 
- \operatorname{ind}(\bD) -\half \int_M \hA(M) \wedge \alpha. \\
\end{eqnarray*}
The integrand $\hA(M) \wedge \alpha$ only contains terms 
of degrees $4p+1$, $p \geq 0$ so the integral over $M$ which 
has dimension $4k+3$ vanishes, and the corollary follows.
\end{proof}

For a self-adjoint first order elliptic operator D the reduced eta
invariant is defined by 
$$
\bareta_D = \eta_D \mod 1 \in \Re/\Za.
$$
We use the same notation as above for the reduced eta invariant of 
twisted Dirac operators, and the for their difference functions with 
respect to variations of the spin structure. In terms of $\bareta$ we
can formulate the following corollary.

\begin{cor}
Under the conditions in Theorem \ref{eta1} we have 
\begin{enumerate}

\item[(a)] $\Delta \bareta(\sigma,\delta;E) = 
-\half \int_M \hA(M) \wedge \chern(\n^E) \wedge \alpha \in \half \Za$, 

\item[(b)] $\Delta^2 \bareta(\sigma,\delta',\delta;E) = 0$.

\end{enumerate}
\end{cor}

\section{The Rokhlin function}

We are now going to prove the theorem on the Rokhlin 
invariant stated in the introduction. We begin by 
stating in more detail the formula by Lee and Miller expressing
the Rokhlin function in terms of eta invariants. 
\begin{thm}\cite{Miller/Lee}, \cite{Liu}, \cite{Liu/Zhang}
\begin{equation} \label{roketa2}
R_M(\sigma) = - \eta_{\operatorname{Hirz}} +
8 \sum_i b_i \eta(\sigma; Z_i) \mod 16,
\end{equation}
where the $b_i$ are integers and the $Z_i$ are tensor bundles, 
that is the $Z_i$ are bundles associated to the orthonormal 
frame bundle $SO(M)$ through representations $\rho_i$ of $SO(n)$.
\end{thm}
This theorem is based on the following formula expressing the
$L$-genus in terms of twisted $\hA$-genera.
\begin{equation} \label{genera}
L(M) = 8 \sum_i b_i \hA(M) \chern(\n^{Z_i}).
\end{equation}

Using Proposition \ref{alphaprop} we see that Theorem \ref{main} 
in the introduction is equivalent to the following.
\begin{thm}
If $\delta \in H^1(M;\Za_2)_0$ then 
$\Delta R_M(\sigma,\delta) = 0 \mod 8$.
\end{thm}

\begin{proof}
From (\ref{roketa2}) we have 
$$
\Delta R_M(\sigma,\delta) = 
8 \sum_i b_i \Delta\eta(\sigma, \delta; Z_i) \mod 16.
$$
Equations (\ref{APS2}) and (\ref{APSintegral}) tells us that 
$$
\Delta R_M(\sigma,\delta) = 
-8 \sum_i b_i \operatorname{ind}(\bD^{Z_i})
-\half \sum_i 8 b_i \int_M \hA(M) \wedge \chern(\n^{Z_i}) \wedge \alpha
\mod 16
$$
and using (\ref{genera}) we see that
$$
\Delta R_M(\sigma,\delta) = 
-8 \sum_i b_i \operatorname{ind}(\bD^{Z_i})
-\half \int_M L(M) \wedge \alpha
\mod 16.
$$
Since $M$ has dimension $8k+3$ and $L(M) \wedge \alpha$ only contains 
terms of degrees $4p+1$, $p \geq 0$, the integral vanishes and we are
left with 
$$
\Delta R_M(\sigma,\delta) = 
-8 \sum_i b_i \operatorname{ind}(\bD^{Z_i}) \mod 16
$$
which proves the theorem.
\end{proof} 

From the relations (\ref{dim3}) and (\ref{dim11}) we get the 
following corollary.

\begin{cor}
If $\delta \in H^1(M;\Za_2)_0$ and $\dim M = 3$ then 
$$
\Delta R_M(\sigma,\delta) = -8 \operatorname{ind}(\bD) \mod 16
$$
and if $\dim M = 11$ then
$$
\Delta R_M(\sigma,\delta) = 8 \operatorname{ind}(\bD^{TM}) \mod 16
$$
where $\bD$ is the operator on $M \times I$ introduced in the proof 
of Theorem (\ref{eta1}).
\end{cor}

\noindent{\bf Acknowledgements:}
The author wishes to thank Sergey Finashin and Stephan Stolz 
for pointing out errors in an earlier version of this paper. 


\providecommand{\bysame}{\leavevmode\hbox to3em{\hrulefill}\thinspace}

\end{document}